\documentclass[12pt]{article}
\usepackage[utf8]{inputenc}
\usepackage[bottom=2cm,marginparwidth=2.5cm]{geometry}

\usepackage{natbib}
\usepackage[USenglish]{babel}
\usepackage{csquotes}

\usepackage{tabularx}
\newcolumntype{L}[1]{>{\hsize=#1\hsize\raggedright\arraybackslash}X} 

\usepackage{xspace}

\usepackage[hyphens]{url}
\usepackage[colorlinks=false,hidelinks]{hyperref}
\newcommand{\furl}[2]{\footnote{{\sloppy\url{#1}}, accessed #2}}
\usepackage[capitalize]{cleveref}

\usepackage{graphicx}
\usepackage{comment}
\usepackage{endnotes}
\let\footnote=\endnote

\title{Research-Data Management Planning in the German Mathematical Community}
\author{\parbox{\linewidth}{Tobias Boege, Ren\'e Fritze, Christiane G\"orgen\footnote{\url{goergen@math.uni-leipzig.de}}, Jeroen Hanselman, Dorothea Iglezakis, Lars Kastner, Thomas Koprucki, Tabea Krause, Christoph Lehrenfeld, Silvia Polla, Marco Reidelbach, Christian Riedel, Jens Saak, Bj\"orn Schembera, Karsten Tabelow, Marcus Weber}}

\begin{document}

\maketitle

\begin{abstract}
In this paper we discuss the notion of research data for the field of mathematics and report on the status quo of research-data management and planning. A number of decentralized approaches are presented and compared to needs and challenges faced in three use cases from different mathematical subdisciplines. We highlight the importance of tailoring research-data management plans to mathematicians' research processes and discuss their usage all along the data life cycle.
\medskip

\textbf{Keywords}\quad research-data management plan, mathematics, use case, scientific computing, computer algebra, theoretical statistics
\end{abstract}

\section{Introduction}

Scientific progress heavily relies on the reusability of previous results.
This in turn is closely linked to reliability and reproducibility of research, and to the question whether another researcher would arrive at the same result with the same material.  In mathematics proofs, together with references to definitions of mathematical objects and already verified theorems, traditionally contained all the information needed in order to verify results. 
However, the advent of computers has opened up new resources previously deemed impossible, while increasing the need for well-adapted research-data management (RDM). For example, algorithms are now implemented to arrive at new conclusions. The size of examples has exploded several orders in magnitude. And some proofs have become too complicated for even the brightest minds, such that software is consulted for thorough understanding and verification\footnote{\sloppy See e.g.\,the story outlined in \url{https://xenaproject.wordpress.com/2020/12/05/liquid-tensor-experiment}, solved with the lean project \url{https://github.com/leanprover-community/lean-liquid}, both accessed 15/11/2022.}.
Studies \citep{Lejaeghere2016, Schappals2017,stodden:7,Riedel:22} from various fields of applied mathematics show that nowadays many results cannot be easily reproduced and hence verified. 

As we outline in \cref{defs}, there are research data in all subdisciplines of mathematics that need responsible organization and documentation in order to ensure they are handled according to the FAIR principles \citep{FAIR} for sustainable, reproducible, and reusable research. One way to achieve this is via a tailored research-data management plan (RDMP), describing the data life cycle over the course of a project 
\citep{Michener15} and providing guidance to fulfill funding requirements\footnote{\sloppy e.g.\ at the European level \url{https://ec.europa.eu/info/funding-tenders/opportunities/docs/2021-2027/horizon/guidance/programme-guide_horizon_en.pdf}, and at the German national level \url{https://www.dfg.de/download/pdf/foerderung/grundlagen_dfg_foerderung/forschungsdaten/forschungsdaten_checkliste_de.pdf}, both accessed 08/11/2022.}.
In mathematics, it is particularly important to treat the RDMP as a living document \citep[cf.]{Dierkes21} because the mathematical research process is hardly projectable and does usually not follow a standardized collection--analysis--report procedure.  
In subfields with experience in using such documentation, three-fold reports -- at the grant-application stage, as a working document, and as a final report -- have proven useful. We discuss this in \cref{rdm,usecases}, spotlighting examples from different subfields, and conclude this article listing central topics for RDMPs in all areas of mathematics.

\section{Mathematical Research Data}\label{defs}

Following \citet[p.130]{KindlingSchirmbacher2013}, we define \emph{research data} as all digital and analog objects that are generated or handled in the process of doing research\footnote{\sloppy This is in line with the notions employed by the DFG \url{https://www.dfg.de/foerderung/grundlagen\_rahmenbedingungen/forschungsdaten/index.html}, forschungsdaten.info \url{https://www.forschungsdaten.info/themen/informieren-und-planen/was-sind-forschungsdaten}, and the MPG \url{https://rdm.mpdl.mpg.de/introduction/research-data-management} e.g., all accessed 07/11/2022.}. In mathematics, research data thus includes paper publications and proofs therein as well as computational results, code, software, and libraries of classifications of mathematical objects. A non-exhaustive list of possible formats and examples is presented in \cref{tab:researchdata,usecases}. 
The apparent diversity of mathematical research-data formats is also reflected in other characteristics such as their storage size, longevity, and state of standardization \citep{KTTKKI,mrd1stpamm,mrdzbmath,deepFAIR,MaRDIproposal}, leading to RDM needs and challenges which are very specific to the discipline of mathematics.

\begin{table}
    \centering
    \begin{tabularx}{\textwidth}{|L{0.7}|L{1.3}|}
        \hline
        \textbf{Research-data type} & \textbf{Examples of data formats} \\\hline
        Mathematical documents &  PDF, LaTeX, XML, MathML\\\hline
        Notebooks & Jupyter, Mathematica, Pluto\\\hline
        Domain-specific research software packages and libraries & R for statistics, Octave, NumPy/SciPy or Julia for matrix computations, CPLEX, Gurobi, Mosel and SCIP for integer programming, or DUNE, deal.II and Trilinos for numerical simulation\\\hline
        Computer-algebra systems & SageMath, SINGULAR, Macaulay2, GAP, polymake, Pari/GP, Linbox, OSCAR, and their embedded data collections\\\hline
        Programs and scripts & written in the packages and systems above, in systems not developed within the mathematical community, input data for these systems (algorithmic parameters, meshes, mathematical objects stored in some collection, the definition of a deep neural network as a graph in machine learning)\\\hline
        Experimental and simulation data & usually series of states of representative snapshots of an observed system, discretized fields, more generally very large but structured datasets as simulation output or experimental output (simulation input and validation), stored in established data formats (i.e.\,HDF5) or in domain-specific
formats, e.g.\,CT scans in neuroscience, material science or hydrology\\\hline
Formalized mathematics & Coq, HOL, Isabelle, Lean, Mizar, NASA PVS library\\\hline
Collections of mathematical objects & L-Functions and Modular Forms Database (LMFDB), Online Encyclopedia of Integer Sequences (OEIS), Class Group Database, ATLAS of Finite Group Representations, Manifold Atlas, GAP Small Groups Library\\\hline
Descriptions of mathematical models in mathematical modeling languages & Modelica for component-oriented modeling of complex systems, Systems Biology Markup Language (SBML) for computational models of biological processes, SPICE for modeling of electronic circuits and devices
and AIMMS or LINGO as a modeling language for integer programming\\\hline
    \end{tabularx}
    \caption{Mathematical research data comes in a variety of data formats \citep[p.26f]{MaRDIproposal}.}
    \label{tab:researchdata}
\end{table}

One of the most apparent challenges is the question what metadata is sufficient for reusability. We will answer this question partially for the mathematical subfields presented in \cref{usecases}. However, as \citet{deepFAIR} note, \enquote{the meaning and provenance of [mathematical research] data must usually be given in the form of complex mathematical data themselves}. It is thus not surprising that there is no common, standardized metadata format yet. 
A search in the RDA Metadata Standards Catalog\furl{https://rdamsc.bath.ac.uk/subject/Mathematics\%20and\%20Statistics}{07/11/2022} at the time of writing reveals five hits, four from a subfield of statistics and one from economics, none of which could encode information about, say, a computer-algebra experiment. This lack of standardization is in contrast to other disciplines such as the life sciences, where the OBO Foundry\furl{https://obofoundry.org/principles/fp-000-summary.html}{07/11/2022} 
hosts more than one hundred interoperable ontologies to describe and link research results, including common naming conventions \citep{Schober2009, arp2015building}.   

Another important aspect of mathematical research data is its particular data life cycle. Again in contrast, for instance, to the life sciences, where older results can be overruled by new evidence, mathematical results that have been proven true remain true indefinitely. Since they cater for other disciplines such as the physical, social, health or life sciences \citep[Fig.~1 and discussion]{MaRDIproposal}, mathematics has a particular responsibility to science to preserve their results in a sustainable manner. We discuss this aspect and how mathematics can be embedded in interdisciplinary research pipelines in more detail in \cref{ta2,ta4}. \Cref{pure} stresses the role thorough documentation plays in this context, using classifications as an example.
%

\section{Status Quo of RDMPs in Mathematics}\label{rdm}

In the narrower sense of data (rather than research data), it is a common claim in the community at the time of writing that mathematics rarely produces data\footnote{\sloppy Usually, only statistics is mentioned as a data-producing subdiscipline, see e.g.\,\url{https://wissenschaftliche-integritaet.de/kommentare/software-entwicklung-und-umgang-mit-forschungsdaten-in-der-mathematik}, accessed 07/11/2022.} and that the few data available need no particular management\footnote{See e.g.\,the inofficial document \url{https://www.math.harvard.edu/media/DataManagement.pdf}, accessed 07/11/2022.}. This is often based on an interpretation of data being something computational, and mathematics being a discipline which is very much paper rather than computer based. Anecdotal evidence suggests that this view is widely established, that there is little knowledge about general RDM, that existing local facilities are hardly used, and that RDMPs are not a standard tool at any stage of the research process. \citet{MaRDIproposal} has identified the need to build common infrastructures for all subdisciplines of mathematics, and mathematics-specific DFG guidelines for FAIR research data will be developed in the foreseeable future\furl{https://www.dfg.de/foerderung/info_wissenschaft/2022/info_wissenschaft_22_25/index.html}{07/11/2022}.

Now, the question of what these guidelines should be is not trivial. 
A large number of questions from a general RDMP catalogue\footnote{For instance the current questionnaire supplied by the DFG-funded research-data management organiser RDMO \url{https://github.com/rdmorganiser/rdmo-catalog/releases/tag/1.1.0-rdmo-1.6.0}, accessed 07/11/2022}, are
irrelevant for a community which produces foremostly theoretical results. For instance, for mathematicians the cost of producing data is rarely relevant -- unlike e.g.\ in the life sciences where data might have to be collected in the field. In the same vein, ethical or data-protection questions most often do not play a role, safe for, for instance, industry collaborations or studies conducted in didactics. Large parts of the community have little training in legal aspects as for example formulae cannot be assigned proprietary rights. In order to avoid the impression that thus all general RDMP questions apply only to sciences different than mathematics, it is imperative to design bespoke catalogues of questions. These should a) use unambiguous language, for instance using the term \enquote{research data} rather than the more specific \enquote{data} which many mathematicians do not handle in their research, and b) avoid superfluous topics while at the same time including sufficient detail, for instance, for mathematics' metadata- and preservation needs identified in the previous section. Now, rather than endeavoring to find a one-size-fits-all solution, in the subsequent section we identify important RDM questions for a number of use cases which are known to the authors -- focusing on metadata, software, data formats and size, versioning, and storage -- and provide those with what we consider to be sensible answers.

We use the remainder of this section to report on two RDM solutions implemented in DFG-funded Collaborative Research Centres (CRC).

The CRC 1456 \enquote{Mathematics of Experiment} includes 17 scientific projects in applied mathematics, computer science, and natural sciences such as biophysics and astronomy, aiming to improve the analysis of experimental data. The research data here are extremely diverse (e.g.\ mathematical documents, notebooks, programs, simulation data or experimental measurements) and their handling is supported by the CRC's dedicated infrastructure project. 
In regular RDM meetings, four themes are recurrent. First, reusage scenarios: especially in interdisciplinary research the same datasets may be processed or used by different groups; documentation, curation, and publication should be tailored to those groups' needs. Second, reproducibility, both computationally and practically in data recreation. Third, metadata: finding accurate descriptors to help the user understand cross-scientific research data. And fourth, visibility: receiving recognition for stand-alone research data beyond a journal publication is hard.
This last topic is usally not part of a standard set of RDMP questions but aims to provide an incentive to increase the effort in research-data creation, publication, and curation.

The CRC 1294 \enquote{Data Assimilation} includes 15 interdisciplinary research projects focusing on the development and integration of algorithms e.g.\,in earthquake prediction, medication dosing, or cell-shape dynamics. 
Researchers are thus confronted both with diverse research data and varying cultural data-handling habits. A central project supports their RDM, and IT infrastructure to facilitate collaborative work and knowledge perpetuation to advance good scientific practice are provided. In particular, the CRC designed an RDMP template in collaboration with the University of Potsdam's research-data group. This covers policies and guidelines, legal and ethical considerations, documentation, and dataset-specific aspects. A vital component of the training is then the classification of the digital objects that are reused and created by the individual researchers. This helps them to develop tailored strategies to improve the quality and reproducibility of published results and to sensitize their research-data handling throughout the data life cycle.

\section{Use Cases}\label{usecases}

We now consider four very different mathematical use cases and discuss their particular research-data needs. Central in these expositions for us is to find out how, using RDMPs, we can provide the best, case-specific guidance to make a project reusable.

\subsection{Applied and Interdisciplinary Mathematics}\label{ta4}

In numerous scientific fields real-world problems are simplified, e.g.\ to experiments, and subsequently described in abstract ways using mathematical models. If a model is combined with input data, it forms a concrete instance of such a problem. With the help of algorithms, the input data is then transformed into output data. Following validations, the interpretation of outputs provides the solution of the initial problem in a so-called Model-Simulation-Optimization workflow \citep[p.77]{MaRDIproposal}. For complete RDM, such workflows should be documented in detail as part of an RDMP.

A standard RDMP questionnaire includes some guidance for the documentation of workflows, such as the main research question, involved disciplines, tools, software, technologies, processes, research-data aspects, and reproducibility. Using this as a template, a tailored questionnaire is currently being developed within the framework of MaRDI\furl{https://www.mardi4nfdi.de}{22/11/2022} to document workflows in detail. This is divided into four sections dealing with the problem statement (object of research, data streams), the model (discretization, variables), the process information (process steps, applied methods), and reproducibility. It is aimed at all disciplines and differs only slightly in whether a theoretical or experimental workflow is documented. The central element of the questionnaire is to establish connections between different steps of the research process in order to improve interoperability of research data. The description of an individual process step, for example, requires the assignment of the relevant input and output data, the method and the (software) environment. At the same time, the documentation of the methods, software, input and output data requires persistent identifiers (e.g.\,Wikidata, swMATH, DOI) in addition to topic-dependent information.

We consider the documentation of a concrete workflow combining archaeology and mathematics as an example. This is based on \citet{Kostre22}, was created by Margarita Kostre independently afterwards, and described in personal communication as \enquote{very helpful for the reflection of the own work}. The author commented that she will use workflow documentation in the future again, as she believes it facilitates interdisciplinary communication, e.g.\ about the status of a project, its goal, and data transfer, it provides better clarity in larger collaborations and allows colleagues to enter a project more easily. The aim of this work is to understand the Romanization of Northern Africa using a susceptible infectious epidemic model. On the process level, the workflow starts with data preparation, e.g.\,collecting, discretizing, and reducing archaelogical data. Once a suitable epidemics model is found, the inverse problem is solved to determine contact networks and spreading-rate functions. Subsequent analysis allows the identification of three different possibilities of the Romanization of Northern Africa. The detailed documentation can be found on the MaRDI Portal\furl{https://portal.mardi4nfdi.de/wiki/Romanization_spreading_on_historical_interregional_networks_in_Northern_Tunisia}{08/11/2022}.

\subsection{Scientific Computing}\label{ta2}
While research in pure mathematics strives to determine an ultimate truth, applied or computational mathematics in majority need to deal with approximations to reality: models are usually expressed in terms of real or complex numbers and only finite subsets of these can actually be implemented on computer hardware. 
Consequently, the result of a computation depends on the format of the finite-precision numbers used and on the specific hardware executing the computations, making a detailed documentation of the computer-based experiment crucial and reusability of code a must-have \citep{morFehHHetal16}. 
Thus, the input data and results of a computer experiment and also the precise implementation (code, software, and hardware) of the algorithms used are important research data.

Absence of such details in documentation makes applied mathematics face the same reproducibility issues \citep[e.g.]{bangerth14} as other scientific fields. Still mathematical algorithms make up the foundation of many computational experiments, for instance as solvers for linear systems of equations, eigenvalue problems, or optimization problems, and are thus at the heart of science today. This responsibility calls for rigorous RDM and documentation in RDMPs.

The main difficulty in establishing RDMPs in scientific computing seems to be in creating incentives to adhere to common standards. In case of a single multi-author paper within a larger project cluster, there are two levels to this question: the funding context and local RDM. Regarding the first, incentives should clearly address reporting requirements and incorporate rewards for sustainable RDM, rather than merely counting publications and citations, to ensure the cluster can {\em stand on the shoulders of giants} instead of {\em building on quicksand}. The beneficiaries here are other researchers in the project and world-wide. Consequently, global RDM needs to answer what is reported where and why. Regarding the local context, for the collaborative work of the authors incentives are far more evident. Thorough RDM, documented in a living RDMP, not only accelerates the paper writing but also improves the reusability of information for future endeavors of the individual authors. Questions center around \enquote{When is the code/data provided? Where in the (local) infrastructure is it stored? By whom? Who is processing it next?}

Consequently, RDMPs should be modularized to enable the single modules to change at their appropriate pace. While the global management rules of a project cluster may not change at all, or at best very slowly, the findings in a single work package may alter the RDMP and thus RDM needs high agility to react to changes. For the software pipeline of an example paper that means: a task-based RDMP, updated as the pipeline evolves, needs to fulfill the requirements of \citet{morFehHHetal16}, while for the project cluster sustainable handover, following \citet[e.g.]{FehHRetal21}, needs to be addressed in the overarching RDMP.

\subsection{Computer Algebra and Theoretical Statistics}\label{pure}
Large parts of the German mathematical community consider themselves as not doing applied work. This includes fields such as geometry, topology, algebra, analysis or number theory, and also mathematical statistics, for instance. However, these researchers increasingly use computers, too, to explore the viability of proof strategies, test their own conjectures or refute established ones. 
As a consequence, \emph{classifications}, the systematic and complete tabulation of all objects with a given property, grow wildly in size and complexity. They give a complete picture for some aspect of a theory and may be used in many ways from the search of (counter)examples over building blocks for constructive proofs to benchmark problems.

For instance, the L-Functions and Modular Forms Database \citep{lmfdb} contains over 4.8TB of data relating objects conjectured to have strong connections by the Langlands program: number fields, elliptic curves, modular forms, L-functions, Galois representations. It includes tens of millions of individual objects and stores the relations between these. Entries contain detailed information on reliability, completeness, and several versions of the code needed to compute them. The database has a public reporting system which allows all users to have visibility of any issues or errors.
Other classification databases targeted at specific audiences are listed at \url{https://mathdb.mathhub.info}.

Computing mathematical objects for classification can often be algorithmically hard and time consuming. But once computed, results are final and independent of the software used. 
With larger computer clusters and better algorithms, it is unreasonable and unsustainable to expect researchers who want to build on existing research to repeat individual computations. This expected reuse increases the need for responsible RDM and triggers challenges which need to be addressed in an RDMP. In particular, four themes are central in this regard. First, how can researchers ensure that their research data is correct and complete? Is the connection of mathematical theory and code sound? Second, how can other researchers access, understand, and reuse the research data? Third, how can one ensure longevity of their research data? And fourth, how can researchers report errors/corrections and upload new versions of research data if necessary?


These questions are neatly addressed in the LMFDB mentioned above. To show how things can go wrong without proper RDM we discuss a classification of all conditional independence structures on up to four discrete random variables, originally published in a series of papers \citep{MatusStudeny,MatusII, MatusIII}. 
Of the $2^{24} = 16\,777\,216$ \emph{a priori} possible patterns of how four random variables can influence each other, only $18\,478$ ($\approx 0.11\%$) are realizable with a probability distribution. \cite{SimecekShortNote} digitized this result and left the field after his PhD in 2007. His research data was deleted in 2021 from his former institute's website -- the only public place which ever held the database\footnote{\sloppy A backup is still available on the Internet Archive at \url{http://web.archive.org/web/20190516145904/http://atrey.karlin.mff.cuni.cz/~simecek/skola/models/}, accessed 14/11/2022.}. It was encoded in a packed binary format which is hard to read, search, and reuse. Some files supporting the correctness of the classification for binary distributions use an unspecified, compiler-specific binary serialization format for floating-point data\footnote{A set of scripts for reading these files is available at \url{https://github.com/taboege/simecek-tools}, accessed 14/11/2022.}. The programs used for the creation and inspection of the database were written in a dialect of the Pascal programming language which has not been maintained since 2006. 
The sparse documentation is in Czech.

This situation can only be fixed by recreating the database from scratch, including proofs. 
An RDMP for this project should emphasize the need to list and document each step of redoing the computations, the use of standard data formats with rich metadata for interoperability and searchability of the database, and ensure future reusability of Šimeček's results.

\section{Discussion and Outlook}
The problem of reusability strongly relates to a phenomenon called \emph{dark data} which \enquote{exists only in the bottom left-hand desk drawer of scientists on some media that is quickly aging} \citep{Heidorn2008}. If research data are not available, they are of course neither traceable nor reusable or FAIR. This phenomenon extends from lost USB sticks and conflicting cloud-based collaboration tools like Dropbox\furl{https://www.dropbox.com}{15/11/2022} and Overleaf\furl{https://www.overleaf.com}{15/11/2022} without local backup to papers containing very condensed complicated proofs that can only be taken up in future work if access to handwritten notes of the authors is also possible. A prime example of this is presented in \cref{pure} where unavailable research data is in stark contrast to the everlasting truth of mathematical results. 
RDMPs are a tool of choice against such issues, serving as a basic measure to organize the full data life cycle.

From the three case studies considered in \cref{usecases}, we derive that RDMPs in mathematics in particular a) stimulate reflection, clarity, and interdisciplinary communication, b) require flexibility and modularization as living RDMPs, and c) facilitate the documentation of iterative computational processes by fostering research-data interoperability and reusability.


We further conclude that archiving and preservation is key in any mathematical subdiscipline. As a very first step to improve the status quo, all research results necessary for reusability (data, code, notes,\ldots) should be stored in a sustainable and findable manner, using resources already documented in an RDMP before a project starts. Ideally in a second step citeable repositories with persistent identifiers for these research data can be chosen, and in a third step these can be annotated with interlinked metadata, implemented via knowledge graphs. Because of the diversity of mathematical research data, the choice of metadata should be made carefully with possible reusage scenarios and interest groups in mind, also documented in an RDMP. If code is part of a publication, thoughts should be given to the detail of documentation and again appropriate citeable long-term repositories. In addition, an RDMP should be used as a tool to identify legal constraints like the compatibility of software licenses before any actual work is conducted. 

\section*{Competing interests}
The authors have no competing interests to declare.

\section*{Author contributions}
All authors made significant contributions to the design of this review as well as drafting and
revising the manuscript. All have approved this final version, agreed to be accountable, and
have approved of the inclusion of those in the list of authors.

\section*{Acknowledgements}

Ren\'e Fritze, Christiane G\"orgen, Jeroen Hanselman, Lars Kastner, Thomas Koprucki, Tabea Krause, Marco Reidelbach, Jens Saak, Bj\"orn Schembera, Karsten Tabelow, and Marcus Weber are supported by MaRDI, funded by the Deutsche Forschungsgemeinschaft (DFG), project number 460135501, NFDI 29/1 \enquote{MaRDI -- Mathematische Forschungsdateninitiative}. Christian Riedel is supported by the DFG, project-ID 318763901 -- SFB1294. Christoph Lehrenfeld is supported by the DFG, project-ID 432680300 -- SFB1456.

The authors are grateful to Margarita Kostre for retrospectively compiling an RDMP for her project \citet{Kostre22} and to Tim Hasler for background and discussion regarding the MATH+ research-data management organizer.

\section*{Author affiliations}
Tobias Boege orcid.org/0000-0001-7284-1827\\
Aalto University \\
René Fritze orcid.org/0000-0002-9548-2238\\
University of Münster\\
Christiane G\"orgen orcid.org/0000-0002-6476-956X\\
Universit\"at Leipzig\\
Jeroen Hanselman orcid.org/0000-0002-1298-0961\\
TU Kaiserslautern \\
Dorothea Iglezakis orcid.org/0000-0002-8524-0569\\
University of Stuttgart\\
Lars Kastner orcid.org/0000-0001-9224-7761\\
Technische Universit\"at Berlin\\
Thomas Koprucki orcid.org/0000-0001-6235-9412\\
WIAS Berlin\\
Tabea H. Krause orcid.org/0000-0001-7275-5830\\
Universit\"at Leipzig\\
Christoph Lehrenfeld orcid.org/0000-0003-0170-8468\\
Georg-August-Universit\"at Göttingen\\
Silvia Polla orcid.org/0000-0002-2395-2448\\
WIAS Berlin\\
Marco Reidelbach orcid.org/0000-0002-1919-1834\\
Zuse Institute Berlin\\
Christian Riedel orcid.org/0000-0001-5154-4153\\
Universit\"at Potsdam\\
Jens Saak orcid.org/0000-0001-5567-9637\\
Max Planck Institute for Dynamics of Complex Technical Systems\\
Björn Schembera orcid.org/0000-0003-2860-6621\\
University of Stuttgart\\
Karsten Tabelow orcid.org/0000-0003-1274-9951\\
WIAS Berlin\\
Marcus Weber orcid.org/0000-0003-3939-410X\\
Zuse Institute Berlin

\theendnotes

\bibliographystyle{agsm2} 
\bibliography{references}
\end{document}